\documentclass[pdflatex,sn-mathphys-num]{sn-jnl}
\usepackage{graphicx}
\usepackage{amsmath,amssymb,amsfonts,amsthm}
\usepackage{xcolor,booktabs,array}
\usepackage{tikz,pgfplots}
\pgfplotsset{compat=1.18}
\newif\ifrevisionblue
\revisionbluefalse
\definecolor{revisionblue}{RGB}{0,65,180}
\newif\ifrevisionred
\revisionredtrue
\definecolor{revisionred}{RGB}{180,0,0}
\AtBeginDocument{\hypersetup{linkcolor=.,citecolor=.,urlcolor=.}}

\newcommand{\R}{\mathbb R}
\newcommand{\C}{\mathbb C}

\newcommand{\J}{\mathbf 1}
\newcommand{\M}{\mathcal M}

\def\ttimes{\,\rotatebox[]{-90}{$\ltimes$}\,}
\newcommand{\tstp}{\mathbin{\ttimes_{\!*}}}
\newcommand{\eye}{\mathsf e}

\newcommand{\GL}{\mathrm{GL}}
\newcommand{\gl}{\mathfrak{gl}}
\DeclareMathOperator{\bcirc}{bcirc}

\DeclareMathOperator{\unfold}{unfold}
\DeclareMathOperator{\fold}{fold}
\DeclareMathOperator{\lcm}{lcm}

\DeclareMathOperator{\rank}{rank}
\DeclareMathOperator{\vct}{vec}
\theoremstyle{thmstyleone}
\newtheorem{theorem}{Theorem}
\newtheorem{prp}[theorem]{Proposition}
\newtheorem{lem}[theorem]{Lemma}
\newtheorem{cor}[theorem]{Corollary}
\theoremstyle{thmstyletwo}
\newtheorem{rem}[theorem]{Remark}
\newtheorem{exa}[theorem]{Example}
\theoremstyle{thmstylethree}
\newtheorem{dfn}{Definition}
\begin{document}
\title[t-Product and t-STP of Cubic Matrices]{t-Product and t-STP of Cubic Matrices with an Application to Hyper-Networked Systems}
\author[1]{\fnm{Daizhan} \sur{Cheng}}\email{dcheng@iss.ac.cn}
\author*[2]{\fnm{Zhengping} \sur{Ji}}\email{zhengping.ji@sorbonne-universite.fr}
\affil[1]{\orgdiv{Key Laboratory of Systems and Control}, \orgname{Academy of Mathematics and Systems Science, Chinese Academy of Sciences}, \orgaddress{\city{Beijing}, \postcode{100190}, \country{China}}}
\affil*[2]{\orgdiv{Laboratoire Jacques-Louis Lions}, \orgname{Sorbonne Universit\'{e}}, \orgaddress{\city{Paris}, \postcode{75005}, \country{France}}}
\abstract{Control systems with tensor-valued state transitions require a product that specifies both how coefficients act on each frontal slice and how different slices interact. This paper develops a t-semi-tensor product (t-STP) on cubic matrices that retains the circular coupling of the t-product while allowing rectangular coefficient slices to act on a fixed state space. The construction combines the dimension-keeping semi-tensor product (DK-STP) bridge with circular convolution, overcoming the absence of cross-slice coupling in a slice-wise DK-STP. It provides a compact coefficient description of a structured class of dynamical operators, with fewer stored entries when the coefficient slices have fewer columns than rows. For a fixed number of frontal slices, we establish associative algebra and module structures and describe the associated Lie algebra and Lie groups. These structures make coefficient composition and exponential evolution consistent, while equivalent classical matrix realizations connect the tensor formulation to control analysis of cubic matrix-based dynamics. A specified supply-network game illustrates how the construction organizes interacting chain flows, reproduces the classical trajectories, and supports a globally convergent payoff-gradient adjustment law with explicit damping. The example quantifies coefficient storage while clarifying that the state dimension is unchanged and that the same economy is available to a classical implementation retaining the factorization.}
\maketitle
\clearpage

\section{Introduction}\label{S1}

With recent advances in data science and computational methods, the analysis of higher-order data structures, particularly multiway arrays or hypermatrices, has attracted considerable attention because of their applications in signal processing \cite{la}, computer vision \cite{va}, and neuroscience \cite{bec}. Once a basis is fixed, a tensor can be represented by its components in a hypermatrix \cite{boo86}. A real hypermatrix of order $d$ is an array $A=(a_{i_1,\ldots,i_d})$, where $i_k\in\{1,\ldots,n_k\}$ for $k=1,\ldots,d$. Equivalently, it is a map from the product of these index sets to $\R$. We denote the space of such arrays by $\R^{n_1\times\cdots\times n_d}$ and refer to \cite{lim13} for a general introduction.

Among higher-order arrays, those of order three provide a direct extension of matrices and occur frequently in data analysis \cite{kil11}. In this paper, we call them \emph{cubic matrices}; their three dimensions need not be equal. We consider arrays
\begin{equation}\label{1.2}
 A=(a_{i,j,k}),\qquad i\in\{1,\ldots,m\},\quad j\in\{1,\ldots,n\},\quad k\in\{1,\ldots,s\},
\end{equation}
where $m,n,s$ are positive integers. The space $\R^{m\times n\times s}$ is the main setting of our study.

Several product operators can be defined on cubic matrices,  and by iteratively applying the product, one can derive dynamics on cubic matrix sets. The choice of product determines which interactions the resulting algebra and dynamical model can represent. For control applications when the state transition law is stored in data of multiway arrays, it is useful to have a product that admits a clear matrix representation, supports repeated multiplication, and accommodates the dimensions of the coefficient and state arrays. These considerations set the basis for comparing the products and their associated dynamical systems below.

The t-product developed by Kilmer and collaborators \cite{kil11,kil13} provides an effective framework for operations on cubic matrices. It represents tensor multiplication through a block-circulant matrix and permits structured computations while retaining the organization of the data. Related t-product methods have been used in imaging and tensor approximation \cite{arz18,kil13}, and tensor functions have also been developed within this framework \cite{lun20}. Kernfeld et al.\ \cite{ker15} generalize tensor products by replacing the Fourier transform with an arbitrary invertible linear transform. Within such tensor-product algebras, truncated tensor singular value decompositions also provide optimal approximation results and economical representations of multiway data \cite{kil21}. These developments show why the choice of multiplication matters for both analysis and representation.

In particular, t-product-based dynamical and control systems have been studied as structured models with tensor-valued states \cite{hoo21,mao24}. Their matrix representations connect tensor operations with familiar tools from linear systems theory. For order-three tensors, a dynamical system over cubic matrices collects several sets of state variables in a single array. Each frontal rectangular slice may represent a subsystem, while multiplication describes interactions among these subsystems without discarding their organization in the state representation. For example, different slices may collect real-valued variables associated with different markets or interacting parts of a network. The product used in the model then specifies how the slices influence one another.

Tensor methods and hypermatrix constructions have been applied to biological systems \cite{yah19}, image processing \cite{ben17}, and control \cite{cche21b,mao24}. These applications involve different tensor products and modeling assumptions; for example, the multilinear systems in \cite{cche21b} use the Einstein product. A broader introduction to tensor-based dynamical systems and the products used in their construction can be found in \cite{cche24}.

Despite these advantages, the t-product requires matching inner dimensions and equal numbers of frontal slices \cite{kil13}. In particular, a left-acting autonomous t-product system uses square frontal slices in its coefficient tensor, although its state slices may be rectangular. The question addressed here is how to retain this slice coupling while describing the state operator by rectangular coefficient slices. In this paper, we propose a construction that separates two modeling choices: a dimension-dependent bridge handles the inner dimensions, and circular convolution specifies the coupling between slices. Whereas the transform-based generalization in \cite{ker15} changes the transform used to define multiplication, our construction retains circular coupling and modifies the multiplication between coefficient and state slices. We study the resulting coefficient factorization and the algebraic conditions that make its use consistent in products, trajectories, and compositions. Comparing the constructions also clarifies which interaction patterns each can represent.

To generalize the t-product and relax the inner-dimension requirement, we first extend the dimension-keeping semi-tensor product (DK-STP) of matrices to cubic matrices and examine the resulting control systems. The DK-STP is a variation of the classical semi-tensor product (STP) of matrices. The STP generalizes conventional matrix multiplication to factors with arbitrary inner dimensions \cite{che12}, with engineering applications surveyed in \cite{li18}, and has been used to analyze Boolean and finite-valued networks \cite{che11,ji}, finite games \cite{che21,che22,ji22}, and dimension-varying systems \cite{che23}. When dimensions do not match, however, the classical STP can enlarge the resulting arrays. The DK-STP was introduced to retain the outer dimensions of the product \cite{che24}. For matrices $A$ and $B$, it has the form $A\Psi B$, where the bridge matrix $\Psi$ is determined by their inner dimensions. This provides a prescribed way to combine factors whose dimensions are incompatible under ordinary multiplication. Bridge matrices also arise in observer-based realizations of control systems \cite{che26}; here their role is to prescribe a coefficient action while keeping the chosen tensor state space fixed.

One motivation for this DK-STP-based construction is to reduce the number of coefficients used to describe a structured dynamical operator. For a state $X\in\R^{m\times q\times s}$, a t-product coefficient tensor has $m^2s$ entries, whereas a rectangular coefficient tensor $A\in\R^{m\times p\times s}$ has $mps$ entries. When $p<m$, the latter description is smaller. This saving does not reduce the $mqs$ state variables or permit every linear operator on that state space to be represented with fewer independent parameters. Rather, it results from imposing a bridge factorization and a specified coupling structure. The construction is therefore relevant when these restrictions express an acceptable interaction model; in that case, the smaller coefficient array defines the corresponding operator exactly. This differs from choosing Tucker factors \cite{kol09} or a truncated tensor singular value decomposition \cite{kil21} to approximate a given data tensor: the bridge fixes a coordinate-aggregation rule before coefficient fitting. We make this distinction explicit and compare equally structured tensor and classical representations in the numerical example.

The slice-wise DK-STP construction, however, has a limitation for coupled models: each frontal-slice subsystem evolves independently of the other state slices for prescribed inputs. To combine rectangular coefficient slices with interactions between state slices, we introduce the t-semi-tensor product (t-STP), which combines the DK-STP bridge with the circular coupling of the t-product. For a fixed number of frontal slices, this product is associative and reduces to the t-product when the inner dimensions match. It therefore provides an algebraic framework for a class of coupled linear and nonlinear control systems over cubic matrices. At the level of the linear state action, the t-STP is a rectangular coefficient factorization of a structured t-product operator; its purpose is to make this factorization part of a consistent product and dynamical formulation.

Our proposed product also induces algebraic structures on cubic matrices. We study the associated algebra and module structures, the Lie algebra defined by the product commutator, and the Lie group obtained by adjoining an identity and considering invertible elements with scalar coordinate one. Their matrix representations clarify how the tensor constructions relate to conventional algebra and which properties can be transferred between the two settings. The link is the standard passage from associative multiplication to commutators and groups of units \cite{hal03}. Its role here is to make explicit how rectangular tensor coefficients compose, how an identity is supplied when no numerical two-sided identity exists, and how coefficient exponentials act on the state. These results connect the product construction to the evolution operators used in the control systems.

From an application perspective, hypermatrices provide a natural way to organize the multiway interactions of a hypergraph \cite{pe}. Tensor-based hypergraph controllability has also been studied through polynomial dynamics on vertex states \cite{cche21}. Here the state variables are instead assigned to supply chains, so the hypergraph organizes the decisions while the chosen product specifies their interaction \cite{che25}. A supply network involving manufacturers, wholesalers, and markets offers a simple example, with one scalar state associated with each supply chain. We use this setting to construct a  synthetic game in which the 24 chains are the players and their flows evolve by payoff-gradient adjustment with explicit damping. This example shows how the proposed coefficient description can be used in a concrete coupled model: we specify the payoffs, compare tensor and classical trajectories and coefficient storage, and establish convergence to the unique Nash equilibrium. 

The main contributions of this paper are as follows:
\begin{itemize}
\item[(i)] We define a t-STP for cubic matrices with a fixed number of frontal slices, combining a bridge between incompatible inner dimensions with circular coupling between slices. This gives rectangular coefficient arrays a consistent action on a fixed tensor state space, while retaining the cross-slice interactions absent from the slice-wise DK-STP. The construction brings these two features into one associative product and identifies the restricted class of t-product operators that it represents.
\item[(ii)] We establish the corresponding algebra and module structures and describe the associated Lie algebra and Lie group, with explicit attention to the adjoined identity and a faithful matrix representation. These results explain how rectangular coefficients can be composed, exponentiated, and related to invertible evolution operators despite the possible absence of a numerical identity. Their significance is to supply a rigorous link between the proposed product and its dynamical use, through explicit realizations of standard algebraic constructions.
\item[(iii)] We formulate the resulting control systems, derive their equivalent classical realizations, and give a numerical supply-network game illustrating the coefficient-storage advantage and the restrictions imposed by the representation. The realizations make standard matrix analysis available and distinguish a smaller coefficient description from a reduction in state dimension. Explicit coefficients, payoffs, and convergent trajectories demonstrate the construction in a reproducible model; comparison with an equally factorized classical implementation clarifies the scope of the storage benefit.
\end{itemize}

The rest of the paper is organized as follows. Section~\ref{S2} reviews the t-product and the matrix DK-STP, as well as the algebraic properties of the t-product and its block-circulant representation. Section~\ref{S4} recalls t-product-based dynamical and control systems. Section~\ref{S5} extends the DK-STP to cubic matrices and explains the resulting slice-wise decoupling. Section~\ref{S6} introduces the fixed-depth t-STP, establishes its basic properties, and derives the corresponding control systems and classical realizations. Section~\ref{S7} studies the associated Lie algebra and Lie group. Section~\ref{S8} presents the supply-network game and the comparison of representations. Section~\ref{S9} concludes with directions for further research.

We use the following notation. $\M_{m\times n}=\R^{m\times n}$ is the space of real $m\times n$ matrices; $I_n$ is the numerical identity matrix; $\J_n$ is the all-ones column vector. The products $\otimes$, $*$, $\ttimes$, and $\tstp$ denote the Kronecker product, t-product, DK-STP, and t-STP, respectively. Ordinary juxtaposition denotes dimension-compatible matrix multiplication or scalar multiplication. We distinguish the tensor $A$, its unfolding $A_u=\unfold(A)$, and its block-circulant matrix $\Gamma(A)$. The notation $A^{(k)}$ always denotes a frontal slice; product powers use $A^{[k]}$ with an explicit product subscript. The symbol $\eye$ denotes a \emph{formal adjoined identity}, which is not an unspecified rectangular numerical matrix.

\section{Preliminaries}\label{S2}

Two ingredients will be used in our construction of the product of cubic matrices: the t-product describes circular interactions between frontal slices, while the matrix DK-STP handles incompatible inner dimensions. We first recall the t-product and its matrix representation, then introduce the bridge and formal identity needed for rectangular coefficients. This order separates the coupling rule from the dimension rule before combining them in Section~\ref{S6}.

\subsection{t-Product of Cubic Matrices}

For $A\in\R^{m\times n\times s}$, the $k$-th frontal slice is the matrix $A^{(k)}=(a_{i,j,k})_{i,j}\in\M_{m\times n}$. Fixing the first or second index similarly gives a horizontal or lateral slice. We work with frontal slices in a fixed order and define
\begin{equation}\label{cu.2.1}
 A_u=\unfold(A)=\begin{bmatrix}A^{(1)}\\A^{(2)}\\\vdots\\A^{(s)}\end{bmatrix}\in\M_{sm\times n}.
\end{equation}
The inverse map is denoted by $\fold$; its target dimensions are determined by context. Unfolding is a linear bijection and changes only the arrangement of entries. It is distinct from $\Gamma(A)$, which has size $sm\times sn$ and repeats the slices as blocks. A general higher-order array can also be matricized by partitioning its ordered indices into row and column indices and using lexicographic order within each group, but no such additional grouping is needed here.

\begin{dfn}[t-product, \cite{kil13}]\label{dcu.2.1}
Let $A\in\R^{m\times n\times s}$ and $B\in\R^{n\times p\times s}$. Define
\begin{equation}\label{cu.2.2}
\Gamma(A)=\bcirc(A)=
\begin{bmatrix}
 A^{(1)}&A^{(s)}&\cdots&A^{(2)}\\
 A^{(2)}&A^{(1)}&\cdots&A^{(3)}\\
 \vdots&\vdots&\ddots&\vdots\\
 A^{(s)}&A^{(s-1)}&\cdots&A^{(1)}
\end{bmatrix},
\end{equation}
and
\begin{equation}\label{cu.2.4}
 A*B=\fold\bigl(\Gamma(A)\unfold(B)\bigr)\in\R^{m\times p\times s}.
\end{equation}
For $s=1$, $\Gamma(A)=A^{(1)}$. When $m=n=1$, $\Gamma(A)$ is the usual circulant matrix of the vector of entries.
\end{dfn}

\begin{dfn}[Slice-wise transpose]\label{dcu.2.2}
The slice-wise transpose $A^{T_f}\in\R^{n\times m\times s}$ is defined by $(A^{T_f})^{(k)}=(A^{(k)})^{\mathsf T}$. More generally, permuting tensor modes permutes their index ranges along with their positions. We distinguish $T_f$ from the usual t-product transpose, which also reverses the order of slices $2,\ldots,s$ \cite{kil13}; no transpose convention is needed for the algebraic results below.
\end{dfn}

\begin{dfn}\label{d2.3}
Let $\mathcal I_n^s\in\R^{n\times n\times s}$ have first slice $I_n$ and all other slices zero. A tensor $A\in\R^{n\times n\times s}$ is invertible under the t-product if there is $B$ in the same space such that $A*B=B*A=\mathcal I_n^s$.
\end{dfn}
Directly from the definition, $\mathcal I_m^s*A=A=A*\mathcal I_n^s$ for $A\in\R^{m\times n\times s}$. The equivalence between tensor and matrix invertibility is proved in Proposition~\ref{pcu.4.2}.

Let $A\in\R^{m\times n\times s}$. Its Frobenius norm is
\begin{equation}\label{cu.2.6}
 \|A\|_F=\left(\sum_{i=1}^m\sum_{j=1}^n\sum_{k=1}^s a_{i,j,k}^2\right)^{1/2}.
\end{equation}

\subsection{Properties of the t-Product}

We recall the t-product properties needed later and include proofs to specify the conventions. Associativity makes repeated products unambiguous, the module action connects coefficients to states, and the map $\Gamma$ transfers these operations to ordinary matrices. These are the properties used in the dynamical constructions below. Addition and scalar multiplication are componentwise, so $\R^{m\times n\times s}$ is a real vector space of dimension $mns$, and its additive group is Abelian.

For an integer $r$, let $\langle r\rangle_s$ be the unique integer congruent to $r$ modulo $s$ in $\{0,\ldots,s-1\}$. In convolution proofs only, write $\widehat A_r=A^{(r+1)}$ for $r=0,\ldots,s-1$. Definition~\ref{dcu.2.1} gives

\begin{equation}\label{cu.3.7}
 \widehat{(A*B)}_k=\sum_{j=0}^{s-1}\widehat A_{\langle k-j\rangle_s}\widehat B_j.
\end{equation}
This convention uses zero-based indices only to simplify modular sums; the original slice notation remains one-based.

\begin{prp}\label{tcu.3.3}
The t-product is bilinear and associative whenever the inner matrix dimensions match and all factors have the same number of frontal slices.
\end{prp}
\begin{proof}
Bilinearity follows from \eqref{cu.3.7}. For $A\in\R^{m\times n\times s}$, $B\in\R^{n\times p\times s}$, and $C\in\R^{p\times q\times s}$, the $k$-th zero-based slice of $(A*B)*C$ is
\[
 \sum_{\ell=0}^{s-1}\sum_{j=0}^{s-1}
 \widehat A_{\langle k-\ell-j\rangle_s}\widehat B_j\widehat C_\ell.
\]
Expanding $A*(B*C)$ and substituting $j=\langle r-\ell\rangle_s$ gives the same sum. Ordinary matrix associativity completes the proof.
\end{proof}

\begin{cor}\label{ccu.3.4}
The space $V=\R^{n\times n\times s}$ is a unital associative real algebra under the t-product, with identity $\mathcal I_n^s$. Its units form a group. The space $\R^{n\times q\times s}$ is a unital left $V$-module under t-product multiplication.
\end{cor}
\begin{proof}
Proposition~\ref{tcu.3.3} proves associativity and distributivity, and Definition~\ref{d2.3} provides the identity. Products of units are units because $(A*B)^{-1}=B^{-1}*A^{-1}$; inverses and the identity belong to the same space. For the module, bilinearity gives the two distributive laws, associativity gives $A*(B*X)=(A*B)*X$, and $\mathcal I_n^s*X=X$ supplies the unit-action law.
\end{proof}

\begin{prp}\label{pcu.4.1}
The map $\Gamma$ is linear and injective. For $A\in\R^{m\times n\times s}$ and $B\in\R^{n\times p\times s}$,
\begin{equation}\label{cu.4.2}
 \Gamma(A*B)=\Gamma(A)\Gamma(B).
\end{equation}
\end{prp}
\begin{proof}
Linearity follows blockwise, and the first block column recovers $\unfold(A)$, proving injectivity. With zero-based block indices, the $(i,j)$ block of the product on the right is
\[
 \sum_{r=0}^{s-1}\widehat A_{\langle i-r\rangle_s}\widehat B_{\langle r-j\rangle_s}
 =\sum_{\ell=0}^{s-1}\widehat A_{\langle i-j-\ell\rangle_s}\widehat B_\ell
 =\widehat{(A*B)}_{\langle i-j\rangle_s}.
\]
These are exactly the blocks on the left.
\end{proof}

\begin{prp}\label{pcu.4.2}
On $\R^{n\times n\times s}$, $\Gamma$ is a unital algebra homomorphism into $\M_{ns\times ns}$ and an isomorphism onto its block-circulant image. A tensor $A$ is a unit if and only if $\Gamma(A)$ is invertible, in which case
\begin{equation}\label{eq:inverse}
 \Gamma(A^{-1})=\Gamma(A)^{-1}.
\end{equation}
Consequently, its restriction to the group of units is a group homomorphism.
\end{prp}
\begin{proof}
Proposition~\ref{pcu.4.1} proves linearity, injectivity, and multiplicativity, and $\Gamma(\mathcal I_n^s)=I_{ns}$. If $A*B=B*A=\mathcal I_n^s$, their images are inverse matrices. Conversely, the inverse of an invertible matrix is a polynomial in that matrix by the classical Cayley--Hamilton theorem. Hence $\Gamma(A)^{-1}$ belongs to the unital algebra $\Gamma(V)$; write it as $\Gamma(B)$. Multiplicativity and injectivity give $A*B=B*A=\mathcal I_n^s$.
\end{proof}

\begin{rem}\label{rcu.4.3}
Thus $\Gamma$ is a faithful representation, not a map onto the entire matrix algebra. We use its inverse only on its image. This precise statement supplies the algebra, ring, monoid, and group correspondences without requiring a separate universal property.
\end{rem}

\subsection{The DK-STP}
The t-product above couples slices but still requires compatible inner matrix dimensions. This subsection reviews the matrix DK-STP \cite{che24}, which supplies the bridge used to relax that requirement. Bridge matrices also arise in observer-based realizations \cite{che26}; here they are used to prescribe a tensor coefficient action on a fixed state space.

	\begin{dfn}\label{d2.2.1}
		Let $A\in\M_{m\times n}$, $B\in\M_{p\times q}$, and $\ell=\lcm(n,p)$. Their DK-STP is
		\begin{equation}\label{2.2.1}
			A\ttimes B=(A\otimes\J_{\ell/n}^{\mathsf T})(B\otimes\J_{\ell/p})\in\M_{m\times q}.
		\end{equation}
	\end{dfn}
	\begin{prp}\label{p2.2.2}
		The bridge representation is
		\begin{equation}\label{2.2.2}
			A\ttimes B=A\Psi_{n\times p}B,\qquad
			\Psi_{n\times p}=(I_n\otimes\J_{\ell/n}^{\mathsf T})(I_p\otimes\J_{\ell/p}).
		\end{equation}
		The product is bilinear and associative for arbitrary matrix sizes, and agrees with ordinary multiplication when $n=p$.
	\end{prp}
	\begin{proof}
		Factor the first replicated matrix as $A(I_n\otimes\J_{\ell/n}^{\mathsf T})$ and the second as $(I_p\otimes\J_{\ell/p})B$. This gives \eqref{2.2.2}. Bilinearity follows. If $C\in\M_{u\times v}$, both associations of the triple product equal $A\Psi_{n\times p}B\Psi_{q\times u}C$. Finally, $\Psi_{n\times n}=I_n$.
	\end{proof}
	\begin{rem}\label{rem:bridge}
		The bridge is fixed by the two inner dimensions, not fitted separately to the factors. Its unnormalized all-ones convention is used throughout this paper. It specifies how coordinates are aggregated and replicated; choosing this product is therefore a modeling assumption. For example,
		\begin{equation}\label{eq:bridge23}
			\Psi_{2\times3}=\begin{bmatrix}2&1&0\\0&1&2\end{bmatrix}.
		\end{equation}
	\end{rem}
	
{The bridge permits multiplication of rectangular coefficients of a fixed shape, but such an algebra need not contain an identity. An identity will be needed for the constant term of the exponential and for invertible compositions. We therefore state explicitly how to adjoin one.}

	\begin{prp}\label{p2.2.4}
		The fixed-shape space $V=\M_{m\times n}$, with ordinary addition and scalar multiplication and product $A\ttimes B=A\Psi_{n\times m}B$, is a finite-dimensional associative real algebra. If $m\ne n$, it has no numerical two-sided identity. For any associative real algebra $V$, the vector space $\widetilde V=\R\oplus V$ becomes a unital associative algebra under
		\begin{equation}\label{eq:unitization}
			(\alpha,A)(\beta,B)=(\alpha\beta,\alpha B+\beta A+AB),
			\qquad \eye=(1,0),
		\end{equation}
		where $AB$ denotes the product of $V$ in this formula.
	\end{prp}
	\begin{proof}
		The first assertion follows from Proposition~\ref{p2.2.2}. A numerical identity $E$ would satisfy $E\Psi_{n\times m}=I_m$ and $\Psi_{n\times m}E=I_n$, since the corresponding multiplication must fix every matrix of $V$. Rank inequalities imply both $m\le n$ and $n\le m$. For the unitization, distributivity is immediate. The $V$ component of either association of $(\alpha,A)(\beta,B)(\gamma,C)$ is
		\[
		\alpha\beta C+\alpha\gamma B+\beta\gamma A+\alpha BC+\beta AC+\gamma AB+(AB)C,
		\]
		where $(AB)C=A(BC)$. The scalar component is $\alpha\beta\gamma$, and $(1,0)$ is a two-sided identity.
	\end{proof}
	We identify $A\in V$ with $(0,A)$ and write $\alpha\eye+A$ for $(\alpha,A)$. The dimension of $\widetilde V$ is $1+\dim V$, even when $V$ already has a native identity. In that case the formal and native identities remain distinct unless an explicit quotient or identification is stated. Invertibility in $\widetilde V$ means existence of a two-sided inverse in $\widetilde V$. Section~\ref{S7} constructs the associated group and Lie algebra, avoiding any use of an undefined rectangular numerical identity.

\section{Dynamical Systems over Cubic Matrices}\label{S4}

We first use the t-product-based dynamical system as a reference for the two constructions that follow in later sections. An action of a monoid $G$ on a set $V$ is a map $\pi:G\times V\to V$ satisfying $\pi(ab,x)=\pi(a,\pi(b,x))$ and $\pi(e,x)=x$. Such an action is also called an $S$-system. When $V$ is a topological space and each map $x\mapsto\pi(a,x)$ is continuous, it is an action by continuous maps. A continuous-time autonomous system additionally specifies a time-indexed semigroup of such maps. Here the exponential map supplies this semigroup, and inputs give a controlled evolution. This viewpoint connects the algebra of Section~\ref{S2} to trajectories, as powers of the coefficient define its exponential, whose action propagates the initial state.

\subsection{Analytic Functions of Cubic Matrices}

The main function needed for linear dynamics is the exponential. We introduce polynomial and analytic functions through $\Gamma$, following the t-function framework of \cite{lun20}; the accompanying polynomial and spectral identities show precisely which ordinary matrix operations are being transferred.

For $A\in\R^{n\times n\times s}$, define
\begin{equation}\label{eq:tpowers}
 A_*^{[0]}=\mathcal I_n^s,\qquad A_*^{[k]}=\underbrace{A*\cdots*A}_{k\text{ factors}},\quad k\ge1.
\end{equation}
For a scalar polynomial $p(z)=\sum_{k=0}^d a_kz^k$, put $p_*(A)=\sum_{k=0}^d a_kA_*^{[k]}$.

\begin{lem}\label{lem:polynomial}
For every such polynomial, $\Gamma(p_*(A))=p(\Gamma(A))$. In particular, $p(\Gamma(A))=0$ implies $p_*(A)=0$.
\end{lem}
\begin{proof}
Apply Proposition~\ref{pcu.4.2} to each power, including the constant term. The implication follows from injectivity.
\end{proof}
\begin{theorem}[Cayley--Hamilton identity]\label{thm:CH}
If $p$ is the characteristic polynomial of $\Gamma(A)$, then $p_*(A)=0$.
\end{theorem}
\begin{proof}
The matrix Cayley--Hamilton theorem gives $p(\Gamma(A))=0$; apply Lemma~\ref{lem:polynomial}.
\end{proof}

\begin{dfn}\label{dfn:teigen}
A scalar $\lambda\in\C$ is a scalar t-eigenvalue of $A$ if there is a nonzero $x\in\C^{n\times1\times s}$ such that $A*x=\lambda x$. The product and unfolding are extended complex-linearly here.
\end{dfn}
\begin{cor}\label{cor:eigen}
The scalar t-eigenvalues of $A$ are exactly the eigenvalues of $\Gamma(A)$, with eigenvectors related by $x\mapsto\unfold(x)\in\C^{ns}$.
\end{cor}
\begin{proof}
Unfolding the defining equation gives $\Gamma(A)\unfold(x)=\lambda\unfold(x)$. Unfolding is bijective, so the implication also holds in reverse.
\end{proof}
This statement concerns the scalar spectral notion just defined; it is not an identification of all tensor eigenvalue conventions.

\begin{dfn}\label{dfn:analytic}
Let $f$ have a scalar Taylor series at zero with radius of convergence $r>0$:
\begin{equation}\label{eq:taylor}
 f(z)=\sum_{k=0}^{\infty}\frac{f^{(k)}(0)}{k!}z^k,\qquad |z|<r.
\end{equation}
For $A\in\R^{n\times n\times s}$ with $\|\Gamma(A)\|_2<r$, define
\begin{equation}\label{eq:tfunction}
 f_*(A)=\sum_{k=0}^{\infty}\frac{f^{(k)}(0)}{k!}A_*^{[k]}.
\end{equation}
Real Taylor coefficients give a real tensor; otherwise the target is complexified.
\end{dfn}
The matrix series converges absolutely in spectral norm. Its partial sums lie in the finite-dimensional, hence closed, image of $\Gamma$. Therefore \eqref{eq:tfunction} converges and
\begin{equation}\label{eq:functionGamma}
 \Gamma(f_*(A))=f(\Gamma(A)).
\end{equation}
The sufficient condition $\sqrt{s}\|A\|_F<r$ can also be used, since $\|\Gamma(A)\|_2\le\|\Gamma(A)\|_F=\sqrt{s}\|A\|_F$. The condition $\|A\|_F<r$ alone is not sufficient: the two-slice scalar tensor $(0.6,0.6)$ has Frobenius norm less than one, while its circulant matrix has eigenvalue $1.2$.

\begin{exa}\label{exa:functions}
The functions
\[
 \exp_*(A)=\mathcal I_n^s+\sum_{k=1}^{\infty}\frac{A_*^{[k]}}{k!},\qquad
 \sin_*(A)=\sum_{k=0}^{\infty}\frac{(-1)^kA_*^{[2k+1]}}{(2k+1)!}
\]
are defined for every $A$. The series
\[
 \log_*(\mathcal I_n^s+A)=\sum_{k=1}^{\infty}\frac{(-1)^{k+1}}{k}A_*^{[k]}
\]
converges, for example, when $\|\Gamma(A)\|_2<1$. This functional calculus is consistent with the matrix-representation approach to t-functions \cite{lun20}.
\end{exa}

\subsection{t-Product-Based Control Systems}

With the exponential defined, the transition from tensor multiplication to a state equation is direct. The following formulation recalls t-product systems of the type studied in \cite{hoo21,mao24} and makes their classical realization explicit. It provides the baseline for identifying what changes when rectangular coefficient slices are introduced.

Consider the continuous-time system
\begin{equation}\label{5.2.7}
 \begin{cases}
 \dot x(t)=A*x(t)+B*u(t),\\
 y(t)=C*x(t),
 \end{cases}
\end{equation}
where $x(t)\in\R^{n\times q\times s}$, $u(t)\in\R^{p\times q\times s}$, $A\in\R^{n\times n\times s}$, $B\in\R^{n\times p\times s}$, and $C\in\R^{r\times n\times s}$. Then $y(t)\in\R^{r\times q\times s}$. Assume the input is locally integrable. The system has a unique absolutely continuous state trajectory on every finite time interval.

\begin{prp}\label{prp:ttrajectory}
The solution with initial state $x(0)$ is
\begin{equation}\label{eq:tsolution}
 x(t)=\exp_*(tA)*x(0)+\int_0^t\exp_*((t-\tau)A)*B*u(\tau)\,d\tau.
\end{equation}
With $X=\unfold(x)$, $U=\unfold(u)$, and $Y=\unfold(y)$, the equivalent matrix system is
\begin{equation}\label{eq:tclassical}
 \dot X=\Gamma(A)X+\Gamma(B)U,\qquad Y=\Gamma(C)X.
\end{equation}
\end{prp}
\begin{proof}
Definition~\ref{dcu.2.1} gives \eqref{eq:tclassical}. The ordinary variation-of-constants formula and \eqref{eq:functionGamma} give \eqref{eq:tsolution} after folding. Unfolding is a linear bijection, proving existence, uniqueness, and equivalence. The integral is taken entrywise.
\end{proof}

\begin{rem}\label{rem:realization}
The state slices in \eqref{5.2.7} may be rectangular. It is the coefficient of its autonomous left action that has square slices. The unfolded state has $nsq$ scalar entries; a complete vectorization yields the state matrix $I_q\otimes\Gamma(A)$.

One may also apply $\Gamma$ to the state, obtaining a larger block-circulant matrix-valued realization. That realization is equivalent only when its initial state and inputs are restricted to the corresponding images of $\Gamma$. The unfolding in \eqref{eq:tclassical} avoids this redundant enlargement. Classical results transfer through this bijection when their admissible input and output structures are preserved. Computational savings require an implementation that exploits the structure, rather than following from a change of notation alone.
\end{rem}

\section{Dynamical Systems over Cubic Matrices via the DK-STP}\label{S5}

We now relax the inner-dimension restriction while keeping the slices separate. Applying the matrix DK-STP to corresponding slices gives the simplest tensor extension and shows what the bridge accomplishes on its own. Its lack of cross-slice state coupling will then motivate the t-STP in Section~\ref{S6}.

\subsection{Kronecker Product and DK-STP of Cubic Matrices}

The Kronecker convention below fixes the ordering of slices, including the explicit repetition lifts discussed later. The slice-wise DK-STP itself is defined first for a common depth.

\begin{dfn}\label{dfn:tkronecker}
For $A\in\R^{m\times n\times s}$ and $B\in\R^{p\times q\times t}$, the tensor Kronecker product $A\otimes B\in\R^{mp\times nq\times st}$ is specified by
\[
 (A\otimes B)_{(i-1)p+i',\,(j-1)q+j',\,(k-1)t+k'}=a_{i,j,k}b_{i',j',k'}.
\]
In particular, its frontal slice numbered $(k-1)t+k'$ is $A^{(k)}\otimes B^{(k')}$.
\end{dfn}

\begin{dfn}\label{dcs.1.3}
For $A\in\R^{m\times n\times s}$ and $B\in\R^{p\times q\times s}$, the slice-wise DK-STP is
\begin{equation}\label{cs.1.3}
 (A\ttimes B)^{(k)}=A^{(k)}\Psi_{n\times p}B^{(k)},\quad k=1,\ldots,s.
\end{equation}
It belongs to $\R^{m\times q\times s}$.
\end{dfn}

\begin{exa}\label{ecs.1.4}
Let
\[
\begin{aligned}
 A^{(1)}&=\begin{bmatrix}1&1\\1&-1\end{bmatrix},&
 A^{(2)}&=\begin{bmatrix}0&1\\1&1\end{bmatrix},\\
 B^{(1)}&=\begin{bmatrix}0&1\\1&0\\-1&1\end{bmatrix},&
 B^{(2)}&=\begin{bmatrix}0&2\\1&-2\\1&0\end{bmatrix}.
\end{aligned}
\]
Using \eqref{eq:bridge23}, the two output slices are
\begin{equation}\label{eq:DKexample}
 (A\ttimes B)^{(1)}=\begin{bmatrix}0&4\\2&0\end{bmatrix},\qquad
 (A\ttimes B)^{(2)}=\begin{bmatrix}3&-2\\4&0\end{bmatrix}.
\end{equation}
Thus a $2\times2$ coefficient slice can act on a $3\times2$ slice without enlarging the outer dimensions of the result.
\end{exa}

\begin{prp}\label{pcs.1.5}
The slice-wise DK-STP is bilinear and associative for any three factors with common depth $s$. For $A\in\R^{m\times n\times s}$, $B\in\R^{p\times q\times t}$, $C\in\R^{n\times r\times s}$, and $D\in\R^{q\times\ell\times t}$,
\begin{equation}\label{cs.1.5}
 (A\otimes B)\ttimes(C\otimes D)=(A\ttimes C)\otimes(B\ttimes D).
\end{equation}
\end{prp}
\begin{proof}
Bilinearity and associativity hold in every slice by Proposition~\ref{p2.2.2}. In \eqref{cs.1.5}, the matrix dimensions match in each relevant slice, so the $(i,j)$ slice on the left is
\[
 (A^{(i)}\otimes B^{(j)})(C^{(i)}\otimes D^{(j)})
 =(A^{(i)}C^{(i)})\otimes(B^{(j)}D^{(j)}).
\]
This is the corresponding slice on the right.
\end{proof}

\begin{rem}[Explicit lifts and their scope]\label{rem:DKlift}
For a multiple $\theta$ of $s$, define the consecutive repetition lift $\mathcal R_s^\theta$ by
\[
 (\mathcal R_s^\theta A)^{(k)}=A^{(1+\lfloor(k-1)/(\theta/s)\rfloor)},\quad k=1,\ldots,\theta.
\]
Thus a depth-three tensor lifts to depth six as $(A^{(1)},A^{(1)},A^{(2)},A^{(2)},A^{(3)},A^{(3)})$. Under Definition~\ref{dfn:tkronecker}, this is $A\otimes\J_{1\times1\times\theta/s}$, with the all-ones factor on the right. One can define a binary slice-wise product of unequal-depth arrays by first lifting both to their least common multiple. The systems and algebraic results below use a fixed depth; in particular, no lifting rule is implicit in a power or differential equation.
\end{rem}

\begin{exa}\label{exa:lift}
For the slice-wise lifted product, let $A$ have three slices
\[
 A^{(1)}=\begin{bmatrix}0&1\\-1&1\end{bmatrix},\quad
 A^{(2)}=\begin{bmatrix}2&1\\1&0\end{bmatrix},\quad
 A^{(3)}=\begin{bmatrix}0&1\\1&0\end{bmatrix},
\]
and let $B^{(1)}=\begin{bmatrix}1&0\\1&0\\-1&1\end{bmatrix}$ and
$B^{(2)}=\begin{bmatrix}1&1\\1&-2\\1&0\end{bmatrix}$. The slices of
$C=(\mathcal R_3^6 A)\ttimes(\mathcal R_2^6 B)$ are
\[
\begin{aligned}
 C^{(1)}=C^{(2)}&=\begin{bmatrix}-1&2\\-4&2\end{bmatrix},&
 C^{(3)}&=\begin{bmatrix}5&2\\3&0\end{bmatrix},\\
 C^{(4)}&=\begin{bmatrix}9&-2\\3&0\end{bmatrix},&
 C^{(5)}=C^{(6)}&=\begin{bmatrix}3&-2\\3&0\end{bmatrix}.
\end{aligned}
\]
This example specifies its lift explicitly and is not used to infer an unequal-depth t-STP semigroup.
\end{exa}

For fixed $m,n,s$, the space $V_D=\R^{m\times n\times s}$ is an associative real algebra under the slice-wise DK-STP. Because multiplication is independent in each slice, invertibility can likewise be considered slice by slice. To describe the corresponding shifted units, let $U=\M_{m\times n}$ with matrix DK-STP and form the product algebra $(\widetilde U)^s$. Write $J=(\eye_U,\ldots,\eye_U)$ for its identity.

\begin{dfn}\label{dcs.2.2}
An element $J+A_0$, $A_0\in V_D$, is invertible if there exists $B_0\in V_D$ such that, in $(\widetilde U)^s$,
\begin{equation}\label{eq:DKunits}
 (J+A_0)(J+B_0)=(J+B_0)(J+A_0)=J.
\end{equation}
Here $J$ is a tuple of formal units, not an ordinary tensor of rectangular identity matrices. Invertibility is asserted for $J+A_0$, not for $A_0$ alone.
\end{dfn}
\begin{prp}\label{pcs.2.3}
The invertible elements of the form $J+A_0$ form a group, the direct product of $s$ copies of the group of units of the form $\eye_U+A^{(k)}$ in $\widetilde U$.
\end{prp}
\begin{proof}
Multiplication in $(\widetilde U)^s$ is componentwise. Products retain scalar coordinate one in each component. The scalar coordinate of an inverse must also be one, so inverses have the stated form. Associativity and the identity come from the product algebra. These observations prove both the group property and the direct-product description.
\end{proof}

\subsection{DK-STP-Based Control Systems}

The coefficient and state tensors can have different column dimensions. The module statement below ensures that a coefficient still acts on the chosen state space and that successive coefficient actions compose consistently. We then write the slice equations to make the resulting coupling pattern visible.

\begin{prp}\label{pcs.2.5}
Let $R_D=(\R^{n\times p\times s},+,\ttimes)$ and $V=\R^{n\times q\times s}$. The action $(A,X)\mapsto A\ttimes X$ makes $V$ a left module over the algebra $R_D$.
\end{prp}
\begin{proof}
The action is bilinear. For each slice, both $A\ttimes(B\ttimes X)$ and $(A\ttimes B)\ttimes X$ equal $A^{(i)}\Psi_{p\times n}B^{(i)}\Psi_{p\times n}X^{(i)}$. Thus the associative action law holds. The action extends to a unital action of the formal unitization by $(\alpha\eye+A)\cdot X=\alpha X+A\ttimes X$.
\end{proof}
Consider
\begin{equation}\label{cs.2.7}
 \dot x=A\ttimes x+\sum_{k=1}^{\mu}B_ku_k(t),\qquad y=C\ttimes x,
\end{equation}
where $x,B_k\in\R^{n\times q\times s}$, $A\in\R^{n\times p\times s}$, $C\in\R^{r\times n\times s}$, and $u_k$ are scalar, locally integrable inputs. For each slice $i=1,\ldots,s$, set $H_i=A^{(i)}\Psi_{p\times n}$. The component form and its solution are
\begin{equation}\label{cs.2.8}
 \dot x^{(i)}=H_i x^{(i)}+\sum_{k=1}^{\mu}B_k^{(i)}u_k(t),\qquad
 y^{(i)}=C^{(i)}x^{(i)},
\end{equation}
\begin{equation}\label{eq:DKsolution}
 x^{(i)}(t)=e^{tH_i}x^{(i)}(0)+
 \sum_{k=1}^{\mu}\int_0^t e^{(t-\tau)H_i}B_k^{(i)}u_k(\tau)\,d\tau.
\end{equation}
The ordinary matrix variation-of-constants formula proves this expression; in particular, the first term has the required initial value.

Nonlinear systems can also be built on the fixed state space $V$:
\begin{equation}\label{cs.2.702}
 \dot x=F_D(x)+\sum_{k=1}^{\mu}G_{D,k}(x)u_k(t),\qquad y=h_D(x).
\end{equation}
For example, use ordinary constant tensors in $V$ and positive-degree DK-STP polynomials in $x$ for $F_D,G_{D,k}$. Convergent power series are allowed on domains where the series and their derivatives converge. These maps take values in $V$, whereas a nonzero scalar constant times a formal unit generally does not. The output map $h_D$ is defined separately with its stated target space. We assume smooth maps and, for local trajectory statements, piecewise continuous inputs.

For this slice-wise polynomial construction, slice $i$ of each vector field depends only on $x^{(i)}$. Thus prescribed inputs give no cross-slice state coupling. An arbitrary smooth map on $V$, or feedback that depends on several slices, need not have this property; the decoupling statement concerns this specific construction. It is appropriate for independent subsystems. We next combine the bridge with circular convolution to introduce cross-slice interactions directly in the product.

\section{The t-STP and t-STP-Based Systems}\label{S6}

The preceding two types of systems provide complementary ingredients: the t-product couples frontal slices, whereas the slice-wise DK-STP permits rectangular coefficient slices. We combine them by inserting the same matrix bridge into every term of the t-product convolution. This gives a common product for the coupled state equation and for repeated coefficient multiplication.

\subsection{The t-STP}

In the unfolding, $I_s\otimes\Psi_{n\times p}$ first maps each slice of the right factor to compatible row coordinates, and $\Gamma(A)$ then mixes the slices by circular convolution. Folding restores the tensor shape without changing its entries.

\begin{dfn}\label{dfn:tstp}
For $A\in\R^{m\times n\times s}$ and $B\in\R^{p\times q\times s}$, define
\begin{equation}\label{eq:tstp}
 A\tstp B=\fold\!\left(\Gamma(A)(I_s\otimes\Psi_{n\times p})\unfold(B)\right)
 \in\R^{m\times q\times s}.
\end{equation}
Equivalently, in the zero-based notation of \eqref{cu.3.7},
\begin{equation}\label{eq:tstpslices}
 \widehat{(A\tstp B)}_k=\sum_{j=0}^{s-1}
 \widehat A_{\langle k-j\rangle_s}\Psi_{n\times p}\widehat B_j.
\end{equation}
\end{dfn}
Since $\Psi_{sn\times sp}=I_s\otimes\Psi_{n\times p}$, \eqref{eq:tstp} can also be written as $\fold(\Gamma(A)\ttimes\unfold(B))$. The equality of bridges follows by substituting $\lcm(sn,sp)=s\lcm(n,p)$ in \eqref{2.2.2}. If $n=p$, this is the t-product. The depth $s$ is fixed throughout the following algebraic constructions.

\begin{exa}\label{exa:tstp}
For the tensors of Example~\ref{ecs.1.4}, the t-STP has slices
\[
\begin{aligned}
 (A\tstp B)^{(1)}&=A^{(1)}\Psi B^{(1)}+A^{(2)}\Psi B^{(2)}
 =\begin{bmatrix}3&2\\6&0\end{bmatrix},\\
 (A\tstp B)^{(2)}&=A^{(2)}\Psi B^{(1)}+A^{(1)}\Psi B^{(2)}
 =\begin{bmatrix}3&2\\-2&8\end{bmatrix},
\end{aligned}
\]
where $\Psi=\Psi_{2\times3}$. The bridge handles the inner-dimension mismatch, and the off-diagonal terms introduce the interaction between slices.
\end{exa}

To use this binary construction in a dynamical system, repeated products must have a well-defined meaning. The next result shows that the bridges remain compatible under reassociation when the depth is fixed. It is the common basis for coefficient powers, module actions, and the Lie commutator studied later.

\begin{prp}[Fixed-depth associativity]\label{pcs.3.4}
For each fixed $s$, the set $\bigcup_{m,n\ge1}\R^{m\times n\times s}$ is a semigroup under the t-STP. The product is bilinear on any two fixed-shape factors.
\end{prp}
\begin{proof}
Closure follows from Definition~\ref{dfn:tstp}. Let $A\in\R^{m\times n\times s}$, $B\in\R^{p\times q\times s}$, and $C\in\R^{u\times v\times s}$. Expanding $((A\tstp B)\tstp C)^{(k+1)}$ gives
\[
 \sum_{\ell=0}^{s-1}\sum_{j=0}^{s-1}
 \widehat A_{\langle k-\ell-j\rangle_s}\Psi_{n\times p}
 \widehat B_j\Psi_{q\times u}\widehat C_\ell.
\]
Expanding $A\tstp(B\tstp C)$ and reindexing as in Proposition~\ref{tcu.3.3} gives the same expression. Matrix associativity proves equality. Bilinearity follows from the same finite-sum formula.
\end{proof}

\begin{rem}[Why depth is fixed]\label{rem:depth}
An unequal-depth binary extension can be defined by lifting the factors with $\mathcal R_s^\theta$ and then applying \eqref{eq:tstp}. It does not give an associative product on all depths. Indeed, for scalar-slice tensors $a=(1)$, $b=(1)$, $c=(1,1)$, this extension, denoted here by $\diamond$, gives
\[
 (a\diamond b)\diamond c=(2,2),\qquad
 a\diamond(b\diamond c)=(4,4).
\]
Both statements follow from unnormalized circular convolution at depth two. Therefore Proposition~\ref{pcs.3.4} does not assert a semigroup over arbitrary third-mode sizes. A different lifting rule would require a separate compatibility and associativity proof.
\end{rem}

\begin{cor}\label{cor:tmodule}
For fixed $n,p,s$, $R_*=\R^{n\times p\times s}$ is an associative real algebra under the t-STP. For each $q$, $V=\R^{n\times q\times s}$ is a left $R_*$-module with action $\pi:R_*\times V\to V$, $\pi(A,X)=A\tstp X$.
\end{cor}
\begin{proof}
The product preserves the shape of $R_*$ and the action preserves the shape of $V$. Bilinearity supplies the additive laws; Proposition~\ref{pcs.3.4} supplies both algebra associativity and the associative action law. As for any such module, the unitized action is $(\alpha\eye+A)\cdot X=\alpha X+A\tstp X$.
\end{proof}

\subsection{t-STP-Based Control Systems}

We now keep the dimensions and scalar input channels of \eqref{cs.2.7} and replace only the slice-wise product by the t-STP. This isolates the effect of introducing circular state coupling. The matrix $H_A$ below then makes both the resulting trajectories and the coefficient restrictions accessible to classical linear-systems analysis:
\begin{equation}\label{cs.4.4}
 \dot x=A\tstp x+\sum_{k=1}^{\mu}B_ku_k(t),\qquad y=C\tstp x.
\end{equation}
Thus $x,B_k\in\R^{n\times q\times s}$, $A\in\R^{n\times p\times s}$, $C\in\R^{r\times n\times s}$, and $y\in\R^{r\times q\times s}$. Let
\begin{equation}\label{eq:HA}
 H_A=\Gamma(A)(I_s\otimes\Psi_{p\times n})\in\M_{ns\times ns}.
\end{equation}

\begin{prp}\label{prp:tstptrajectory}
For locally integrable inputs, the unique trajectory of \eqref{cs.4.4}, with $X=\unfold(x)$, $B_{k,u}=\unfold(B_k)$, and $Y=\unfold(y)$, satisfies
\begin{equation}\label{eq:tstpclassical}
 \dot X=H_AX+\sum_{k=1}^{\mu}B_{k,u}u_k(t),\qquad Y=\Gamma(C)X,
\end{equation}
and
\begin{equation}\label{eq:tstpsolution}
 X(t)=e^{tH_A}X(0)+\sum_{k=1}^{\mu}\int_0^t e^{(t-\tau)H_A}B_{k,u}u_k(\tau)\,d\tau.
\end{equation}
\end{prp}
\begin{proof}
Unfold \eqref{cs.4.4} using Definition~\ref{dfn:tstp}. The output bridge is $\Psi_{n\times n}=I_n$, giving $Y=\Gamma(C)X$. The standard matrix solution gives \eqref{eq:tstpsolution}; folding gives the unique tensor trajectory.
\end{proof}

The matrix solution already determines the trajectory. To express the same propagation intrinsically through the coefficient algebra, form the unitization $\widetilde R_*$. Define $A_{\tstp}^{[0]}=\eye$ and $A_{\tstp}^{[k]}$ by repeated t-STP for $k\ge1$. Choose a submultiplicative norm on this finite-dimensional unital algebra, for example the operator norm of its faithful left regular representation. If $f$ has scalar Taylor radius $r$, the series

\begin{equation}\label{eq:tstpfunction}
 f_{\tstp}(A)=f(0)\eye+\sum_{k=1}^{\infty}\frac{f^{(k)}(0)}{k!}A_{\tstp}^{[k]}
\end{equation}
converges whenever $\|A\|_{\widetilde R_*}<r$. The exponential converges for every $A$ and all finite times. The module action gives
\begin{equation}\label{eq:tstpTensorSolution}
\begin{split}
 x(t)={}&\exp_{\tstp}(tA)\cdot x(0)\\
 &+\sum_{k=1}^{\mu}\int_0^t\exp_{\tstp}((t-\tau)A)\cdot B_k\,u_k(\tau)\,d\tau.
\end{split}
\end{equation}
Indeed, unfolding the action of $A$ is left multiplication by $H_A$; the identity acts as $I_{ns}$, and the convergent power series acts as $e^{tH_A}$. Thus \eqref{eq:tstpTensorSolution} is exactly \eqref{eq:tstpsolution}, with the formal identity acting on the state rather than being mistaken for a state tensor.

\begin{rem}[Advantage and restriction of the factorization]\label{rem:factorization}
Set $K^{(i)}=A^{(i)}\Psi_{p\times n}\in\M_{n\times n}$. Then $H_A=\bcirc(K)$ and $A\tstp x=K*x$. The t-STP therefore gives a rectangular coefficient description of a structured t-product operator. If $p<n$, $A$ has $nps$ entries instead of the $n^2s$ entries of an unrestricted $K$, but each block $K^{(i)}$ has rank at most $p$ and its rows belong to the fixed row space of $\Psi_{p\times n}$. Not every square-slice t-product coefficient has this factorization.

The restriction also affects stability: $\rank H_A\le sp<sn$ when $p<n$, so the autonomous system $\dot X=H_AX$ cannot have an asymptotically stable origin on the full state space. Damping or additional feedback can remove this obstruction; it must be specified separately if it lies outside the bridge-factorized family. Furthermore, $\vct(X)$ evolves under $I_q\otimes H_A$, so the autonomous left action does not mix different columns of $X$.
\end{rem}

\begin{rem}\label{rem:specializations}
When $p=n$, the state operator in \eqref{cs.4.4} reduces to the t-product state operator. Its bank of scalar input channels need not coincide with the tensor-input parameterization of \eqref{5.2.7}; both describe linear input maps once their channels are specified. Replacing the t-STP by the slice-wise DK-STP instead gives \eqref{cs.2.7}, with a different coupling pattern.
\end{rem}

The same product also supplies nonlinear interaction terms: positive-degree powers stay in the state space and mix slices through convolution. For nonlinear systems let $V=\R^{n\times q\times s}$, let $U\subseteq V$ be open, and choose smooth maps $F,G_k:U\to V$ and $h:U\to\R^{r\times\ell\times s}$. Then

\begin{equation}\label{cs.4.7}
 \dot x=F(x)+\sum_{k=1}^{\mu}G_k(x)u_k(t),\qquad y=h(x)
\end{equation}
is a well-defined control-affine system, with local solutions for piecewise continuous inputs. A particular analytic class is obtained from ordinary constant tensors in $V$ and positive-degree t-STP polynomials in $x$, or convergent series on an appropriate open domain. For example, $F(x)=D+a x+b(x\tstp x)$, with $D\in V$ and $a,b\in\R$, takes values in $V$ and introduces cross-slice nonlinear terms. General analytic vector fields need not be scalar-coefficient functions of one algebra element. In particular, the formal constant $f(0)\eye$ in \eqref{eq:tstpfunction} cannot be inserted into \eqref{cs.4.7} as an ordinary rectangular state tensor. The output $h$ is a map, not necessarily a vector field.

\section{Lie Groups and Lie Algebras via t-STP}\label{S7}

The previous section used products and exponentials to describe state evolution. We now examine the algebra of their composition: the commutator records the dependence on the order of coefficient actions, and the group of units contains their invertible exponentials. Fix $m,n,s$ and put $V=\R^{m\times n\times s}$ with the t-STP. The following constructions use the finite-dimensional associative algebra from Proposition~\ref{pcs.3.4}. They extend the rectangular matrix setting of \cite{che24} using the standard relation between associative algebras and matrix Lie structures \cite{hal03}. They concern the product commutator, which should not be identified without further argument with the Lie bracket of nonlinear vector fields in \eqref{cs.4.7}.

\begin{lem}\label{lcs.5.3}
For $A,B\in V$, the block-circulant map satisfies
\begin{equation}\label{cs.5.01}
 \Gamma(A\tstp B)=\Gamma(A)(I_s\otimes\Psi_{n\times m})\Gamma(B)
 =\Gamma(A)\ttimes\Gamma(B).
\end{equation}
Thus it is an injective algebra homomorphism into the rectangular matrix algebra $\M_{sm\times sn}$ equipped with its DK-STP, not in general into a square matrix algebra with ordinary multiplication.
\end{lem}
\begin{proof}
The $(i,j)$ block of the first product on the right is the circular convolution of the slices of $A$ and $B$ with the fixed bridge inserted between them, which is the $(i,j)$ block on the left by \eqref{eq:tstpslices}. The last equality uses $\Psi_{sn\times sm}=I_s\otimes\Psi_{n\times m}$. Injectivity and linearity hold as in Proposition~\ref{pcu.4.1}.
\end{proof}

For the linear coefficient actions $T_A(X)=A\tstp X$ on $\R^{m\times q\times s}$, associativity gives $T_AT_B-T_BT_A=T_{A\tstp B-B\tstp A}$. Thus the following commutator measures the effect of exchanging two such actions and remains representable by a tensor of the same shape.

\begin{prp}\label{pcs.5.1}
The bracket
\begin{equation}\label{cs.5.1}
 [A,B]_{\tstp}=A\tstp B-B\tstp A
\end{equation}
makes $V$ a real Lie algebra, denoted by $\gl_*(m\times n\times s,\R)$. For $m=n$ it is the commutator Lie algebra of the t-product algebra.
\end{prp}
\begin{proof}
Bilinearity and skew-symmetry follow from the product. By associativity, suppressing the product sign in this calculation only, the three terms of the Jacobi identity are
\[
\begin{aligned}
 {}[[A,B],C]&=ABC-BAC-CAB+CBA,\\
 [[B,C],A]&=BCA-CBA-ABC+ACB,\\
 [[C,A],B]&=CAB-ACB-BCA+BAC.
\end{aligned}
\]
Their sum is zero. When $m=n$, the bridge is $I_n$ and the t-STP equals the t-product.
\end{proof}

To pass from infinitesimal coefficients to invertible compositions, the zero coefficient must correspond to an identity action. Adjoin one formal identity by $\widetilde V=\R\oplus V$ with multiplication \eqref{eq:unitization}. Thus
\begin{equation}\label{eq:tunitization}
 (\alpha\eye+A)(\beta\eye+B)
 =\alpha\beta\eye+\alpha B+\beta A+A\tstp B.
\end{equation}
This unitization of the convolution algebra differs from the product of slice unitizations used for $J$ in Definition~\ref{dcs.2.2}. Both identities are formal, but they belong to different algebras.

\begin{dfn}\label{dcs.5.2}
An element $z\in\widetilde V$ is a unit if there exists $w\in\widetilde V$ with $zw=wz=\eye$. Define
\begin{equation}\label{eq:tgroup}
 \GL_*(m\times n\times s,\R)
 =\{\eye+A:A\in V,\ \eye+A\text{ is a unit in }\widetilde V\}.
\end{equation}
\end{dfn}
To study invertibility and smoothness with ordinary matrix tools, we use the left regular representation of the unitized algebra. The additional coordinate records the formal identity, and the first column retains the coefficient itself even when distinct coefficients have the same left action on the original algebra. Let $D=mns$, fix the coordinate map $a(A)=\vct(\unfold(A))\in\R^D$, and let $L_A\in\M_{D\times D}$ satisfy $a(A\tstp B)=L_Aa(B)$. Define
\begin{equation}\label{eq:faithful}
 \rho(A)=\begin{bmatrix}0&0_{1\times D}\\a(A)&L_A\end{bmatrix},\qquad
 \widetilde\rho(\alpha\eye+A)=\alpha I_{D+1}+\rho(A).
\end{equation}

\begin{prp}\label{prp:faithful}
The map $\widetilde\rho$ is a faithful unital algebra representation on $\R^{D+1}$. An element $z\in\widetilde V$ is invertible if and only if $\widetilde\rho(z)$ is invertible as an ordinary matrix. For a specified pair $z,w$, the identity $w=z^{-1}$ is equivalent to $\widetilde\rho(w)=\widetilde\rho(z)^{-1}$. The restriction $\rho$ is an injective Lie algebra homomorphism into $\gl(D+1,\R)$.
\end{prp}
\begin{proof}
Associativity gives $L_{A\tstp B}=L_AL_B$. Multiplying the two block matrices yields
\[
 \rho(A)\rho(B)=
 \begin{bmatrix}0&0\\L_Aa(B)&L_AL_B\end{bmatrix}
 =\rho(A\tstp B).
\]
The first column recovers $a(A)$, so $\rho$ is injective; the upper-left entry of $\widetilde\rho(\alpha\eye+A)$ recovers $\alpha$. Formula~\eqref{eq:tunitization} now proves unital multiplicativity. In fact $\widetilde\rho(z)$ is the matrix of left multiplication by $z$ on $\widetilde V$ in the basis beginning with $\eye$. If this matrix is invertible, there is a unique $w$ with $zw=\eye$. Then $z(wz-\eye)=0$, and injectivity of left multiplication by $z$ implies $wz=\eye$. The reverse implication follows by representing an inverse. The claim about a specified pair follows by faithfulness. Finally, multiplicativity gives $\rho([A,B]_{\tstp})=[\rho(A),\rho(B)]$.
\end{proof}

\begin{theorem}\label{tcs.5.6}
The set \eqref{eq:tgroup} is a Lie group of dimension $D=mns$. Its Lie algebra is $\gl_*(m\times n\times s,\R)$. Moreover,
\begin{equation}\label{eq:groupembedding}
 \eye+A\longmapsto I_{D+1}+\rho(A)
\end{equation}
is a faithful smooth group representation and an isomorphism onto an embedded Lie subgroup of $\GL(D+1,\R)$.
\end{theorem}
\begin{proof}
The product of two elements of \eqref{eq:tgroup} has scalar coordinate one by \eqref{eq:tunitization}. The scalar-coordinate map $\alpha\eye+A\mapsto\alpha$ is an algebra homomorphism, so the inverse of an element with scalar coordinate one also has scalar coordinate one. The identity is $\eye$, and associativity is inherited from $\widetilde V$. This proves the group axioms.

By Proposition~\ref{prp:faithful}, the group is the open subset
\[
 \{\eye+A:\det(I_{D+1}+\rho(A))\ne0\}
\]
of the affine space $\eye+V$. Multiplication is polynomial. Matrix inversion is smooth on the invertible matrices, and the lower part of the first column of the inverse recovers the coordinate vector of its $V$ component; inversion in this group is therefore smooth. Its tangent space at $\eye$ identifies with $V$. The differential of \eqref{eq:groupembedding} is $\rho$, so the induced Lie bracket is the commutator \eqref{cs.5.1}, by Proposition~\ref{prp:faithful}.

The representation is injective and has a smooth inverse on its image, obtained from the same first column. Its image is $(I_{D+1}+\rho(V))\cap\GL(D+1,\R)$, an open subset of an affine subspace with the inherited embedded manifold structure. Thus it is an embedded Lie subgroup and the asserted isomorphism follows.
\end{proof}

This group contains the evolution operators already used in Section~\ref{S6}: for $A\in V$, the element $\exp_{\tstp}(tA)$ has scalar coordinate one and inverse $\exp_{\tstp}(-tA)$. Its image under $\widetilde\rho$ is $\exp(t\rho(A))$, and its module action gives the homogeneous solution in \eqref{eq:tstpTensorSolution}. Thus the group and the dynamical formulation describe the same coefficient exponential in two representations.

\begin{rem}\label{rem:representationscope}
The representation \eqref{eq:faithful} proves existence of a faithful classical realization of the algebraic structures. Its size $D+1$ is not a computational compression claim. The smaller state-action representation $A\mapsto H_A$ from \eqref{eq:HA} is useful for dynamics but need not be faithful for arbitrary rectangular shapes. These two purposes must be distinguished. In particular, no normalized squaring map is needed to establish Theorem~\ref{tcs.5.6}.
\end{rem}

\section{Application to Hyper-Networked Evolutionary Games}\label{S8}

We apply the preceding construction to a specified continuous-time payoff-gradient game on a supply hypernetwork. The synthetic example illustrates the utility of the operator structure proposed in this paper by following three steps: organize the chain decisions as a tensor, specify an interaction and payoffs that yield the state equation, and compare its tensor and classical realizations. Hypergraph signal processing \cite{pe} and hypergraph controllability \cite{cche21} provide related uses of tensor representations; here the focus is a bridge-factorized interaction between chain flows.

\subsection{Network, players, and state}

The hypergraph specifies which three actors form a chain; its dual allows those chains to be indexed as decision units. This separates the network incidence from the numerical law governing the flow on each chain.

Consider two manufacturers $x_1,x_2$, three wholesalers $y_1,y_2,y_3$, and four markets $z_1,\ldots,z_4$. The hypergraph $H$ has nine vertices and 24 hyperedges
\begin{equation}\label{eq:hyperedges}
 E_{ijk}=\{x_i,y_j,z_k\},\qquad
 i=1,2,\quad j=1,2,3,\quad k=1,\ldots,4.
\end{equation}
Every manufacturer--wholesaler--market path in Figure~\ref{fig:supply} represents one hyperedge. The dual hypergraph $H^*$ has one vertex for each chain. To make its labeling explicit, let
\begin{equation}\label{eq:chainindex}
 a(i,j,k)=j+3(k-1)+12(i-1).
\end{equation}
Its nine hyperedges, corresponding to the original actors, are
\[
\begin{aligned}
 e^x_i&=\{a(i,j,k):j=1,2,3;\ k=1,\ldots,4\},\\
 e^y_j&=\{a(i,j,k):i=1,2;\ k=1,\ldots,4\},\\
 e^z_k&=\{a(i,j,k):i=1,2;\ j=1,2,3\}.
\end{aligned}
\]
These formulas fix all incidences without a second, ambiguous vertex ordering.

\begin{figure}[t]
\centering
\begin{tikzpicture}[x=1.35cm,y=1.05cm,every node/.style={font=\small}]
\node (x1) at (1,2.5) {$x_1$};\node (x2) at (3,2.5) {$x_2$};
\node (y1) at (0.5,1.25) {$y_1$};\node (y2) at (2,1.25) {$y_2$};\node (y3) at (3.5,1.25) {$y_3$};
\node (z1) at (0,0) {$z_1$};\node (z2) at (1.33,0) {$z_2$};\node (z3) at (2.67,0) {$z_3$};\node (z4) at (4,0) {$z_4$};
\foreach \i in {1,2}{\foreach \j in {1,2,3}{\draw[gray!55] (x\i)--(y\j);}}
\foreach \j in {1,2,3}{\foreach \k in {1,2,3,4}{\draw[gray!55] (y\j)--(z\k);}}
\draw[very thick] (x1)--(y2)--(z3);
\node[anchor=west] at (4.25,2.5) {Manufacturers};
\node[anchor=west] at (4.25,1.25) {Wholesalers};
\node[anchor=west] at (4.25,0) {Markets};
\end{tikzpicture}
\caption{Supply-network incidence schematic. Each two-link path specifies a three-actor hyperedge; the thick path highlights $E_{123}$. The nine actors determine the network, while its 24 chains are the players of the illustrative game.}\label{fig:supply}
\end{figure}

Each chain is treated as a coordinated decision unit and chooses a scalar flow $w_{ijk}\ge0$. We collect these 24 decisions in a \emph{real-valued} cubic matrix $W\in\R^{3\times2\times4}$:
\begin{equation}\label{eq:flowstate}
 W^{(k)}=\begin{bmatrix}
 w_{11k}&w_{21k}\\w_{12k}&w_{22k}\\w_{13k}&w_{23k}
 \end{bmatrix},\qquad X=\unfold(W)\in\R^{12\times2}.
\end{equation}
The vector $v=\vct(X)$ is column-major, with $v_{a(i,j,k)}=w_{ijk}$. The state therefore has 24 scalar coordinates; it is not a tensor of strategy triples.

\subsection{Interaction and payoff functions}

We next choose the interaction coefficients and then construct payoffs whose gradients generate the desired dynamics. This is a continuous-decision example; the finite-game formulations in \cite{che21,che22,ji22} provide related game-control motivation but their convergence results are not assumed here. Choose four numerical $3\times2$ coefficient slices:

\begin{equation}\label{eq:Anumeric}
\begin{aligned}
 A^{(1)}&=\begin{bmatrix}1&0\\1/2&1/2\\0&1\end{bmatrix},&
 A^{(2)}&=\begin{bmatrix}0&1\\1/2&1/2\\1&0\end{bmatrix},\\
 A^{(3)}&=\begin{bmatrix}1/2&1/2\\1&0\\0&1\end{bmatrix},&
 A^{(4)}&=\begin{bmatrix}0&1\\1&0\\1/2&1/2\end{bmatrix}.
\end{aligned}
\end{equation}
With $\Psi=\Psi_{2\times3}$ from \eqref{eq:bridge23}, let $K_i=A^{(i)}\Psi$. Direct multiplication gives
\begin{equation}\label{eq:Knumeric}
\begin{aligned}
 K_1&=\begin{bmatrix}2&1&0\\1&1&1\\0&1&2\end{bmatrix},&
 K_2&=\begin{bmatrix}0&1&2\\1&1&1\\2&1&0\end{bmatrix},\\
 K_3&=\begin{bmatrix}1&1&1\\2&1&0\\0&1&2\end{bmatrix},&
 K_4&=\begin{bmatrix}0&1&2\\2&1&0\\1&1&1\end{bmatrix}.
\end{aligned}
\end{equation}
Consequently,
\begin{equation}\label{cs5.1101}
 \unfold(A\tstp W)=LX,\qquad
 L=\begin{bmatrix}
 K_1&K_4&K_3&K_2\\K_2&K_1&K_4&K_3\\
 K_3&K_2&K_1&K_4\\K_4&K_3&K_2&K_1
 \end{bmatrix}\in\M_{12\times12}.
\end{equation}
The frontal slices represent markets, and the circular coupling in $L$ assumes that the effect depends on their relative index modulo four. This is a deliberate modeling restriction, not a consequence of the hypergraph alone. Different manufacturer columns evolve independently under this left action, although chains within each column interact across wholesalers and markets.

Every row of $L$ sums to 12, so the uncontrolled equation $\dot X=LX$ has the growing mode $e^{12t}\J_{12}$. We specify damping and incentives separately:
\begin{equation}\label{eq:gameTensor}
 \dot W=A\tstp W-13W+\J_{3\times2\times4}.
\end{equation}
The damping $-13W$ need not lie in the same bridge-factorized coefficient family and is not absorbed into $A$. It can also be viewed as a specified feedback term in a system with additive input channels.

To realize \eqref{eq:gameTensor} as a game, define
\begin{equation}\label{eq:N}
 N=I_2\otimes(L-13I_{12})\in\M_{24\times24}.
\end{equation}
Player $a$ chooses $v_a\in[0,\infty)$ and has payoff
\begin{equation}\label{eq:payoffs}
 P_a(v)=v_a+\frac12N_{aa}v_a^2+v_a\sum_{b\ne a}N_{ab}v_b,
 \qquad a=1,\ldots,24.
\end{equation}
All payoff coefficients are explicit through \eqref{eq:Anumeric}--\eqref{eq:N}; no unspecified payoff tensors remain. With adjustment rate one, the payoff-gradient dynamics are
\begin{equation}\label{cs5.12}
 \dot v_a=\frac{\partial P_a}{\partial v_a}
 =1+N_{aa}v_a+\sum_{b\ne a}N_{ab}v_b.
\end{equation}
Thus \eqref{cs5.12} is exactly the vectorization of \eqref{eq:gameTensor}. The own quadratic term penalizes large flow, and the linear interaction terms describe the chosen incentives. This is a continuous-time payoff-gradient adjustment model, not a replicator equation. The payoffs are quadratic, rather than linear, so their gradients give nonconstant feedback.

\begin{prp}\label{prp:gameequilibrium}
The nonnegative orthant is forward invariant under \eqref{cs5.12}. The game \eqref{eq:payoffs} has the unique Nash equilibrium $v^*=\J_{24}$, and every trajectory with nonnegative initial state converges to it. More precisely,
\begin{equation}\label{eq:contraction}
 \|v(t)-\J_{24}\|_\infty\le e^{-t}\|v(0)-\J_{24}\|_\infty,\qquad t\ge0.
\end{equation}
\end{prp}
\begin{proof}
All off-diagonal entries of $N$ are nonnegative, and $N\J_{24}=-\J_{24}$, so $\J_{24}$ is an equilibrium. A matrix with nonnegative off-diagonal entries has a nonnegative exponential: add a sufficiently large scalar multiple of the identity and use the exponential series. Variation of constants then proves nonnegativity for nonnegative initial states and nonnegative constant forcing.

Let $Q=N+I_{24}$. It has nonnegative off-diagonal entries and zero row sums. Thus $e^{tQ}$ is nonnegative and $e^{tQ}\J_{24}=\J_{24}$, so its induced infinity norm is one. Since $e^{tN}=e^{-t}e^{tQ}$, the error $v(t)-\J_{24}=e^{tN}(v(0)-\J_{24})$ satisfies \eqref{eq:contraction}. In particular, $N$ is invertible.

For fixed nonnegative decisions of the other players, $P_a$ is strictly concave in $v_a$ because $N_{aa}\le-11$. Its unique best response is
\[
 \frac{1+\sum_{b\ne a}N_{ab}v_b}{-N_{aa}}>0.
\]
Every Nash equilibrium is therefore interior and solves $Nv+\J_{24}=0$. Invertibility of $N$ gives uniqueness, and $\J_{24}$ satisfies these equations. The contraction estimate proves convergence to this equilibrium.
\end{proof}

\subsection{Numerical realization and coefficient comparison}

The model is now fully specified. We check that its two representations give the same trajectory and then compare the entries needed to store their common interaction operator. These two checks distinguish correctness of the realization from economy of the coefficient description.

Choose
\begin{equation}\label{eq:initial}
 X_{r,1}(0)=\frac{2r-1}{24},\qquad X_{r,2}(0)=\frac{2r}{24},\qquad r=1,\ldots,12.
\end{equation}
With $M=L-13I_{12}$, the exact matrix solution is
\begin{equation}\label{eq:gameSolution}
 X(t)=\J_{12\times2}+e^{tM}(X(0)-\J_{12\times2}).
\end{equation}
It equals the vectorized solution $v(t)=\J_{24}+e^{tN}(v(0)-\J_{24})$. Direct evaluation of the slice formula \eqref{eq:tstpslices} and $LX(0)$ gave a largest entrywise discrepancy of $1.78\times10^{-15}$. At $t=0,0.1,0.5,1,2$, the two trajectory formulas differed entrywise by at most $2.22\times10^{-16}$ in double-precision evaluation. The matrix exponentials were evaluated by scaling to infinity norm at most one, summing the Taylor series to an incremental term norm below $10^{-17}$, and repeated squaring. These are numerical consistency checks, not timing benchmarks. The exact equivalence follows already from Proposition~\ref{prp:tstptrajectory}. Figure~\ref{fig:trajectories} illustrates four of the flows.

\begin{figure}[t]
\centering
\begin{tikzpicture}
\begin{axis}[width=.88\linewidth,height=5.5cm,xlabel={$t$},ylabel={Chain flow},
 xmin=0,xmax=5,ymin=0,ymax=1.04,grid=major,
 legend style={at={(.98,.04)},anchor=south east,font=\small},
 tick label style={font=\small},label style={font=\small},
 every axis plot/.append style={color=revisionblue}]
\addplot[thick,solid] coordinates {
(0.000000,0.041666667) (0.062500,0.331531653) (0.125000,0.472602691) (0.187500,0.548167315) (0.250000,0.594443689) (0.312500,0.627202268) (0.375000,0.653333136) (0.437500,0.675869925) (0.500000,0.696170880) (0.562500,0.714865357) (0.625000,0.732264437) (0.687500,0.748539053) (0.750000,0.763797292) (0.812500,0.778117992) (0.875000,0.791565403) (0.937500,0.804195648) (1.000000,0.816059620) (1.062500,0.827204342) (1.125000,0.837673646) (1.187500,0.847508564) (1.250000,0.856747580) (1.312500,0.865426816) (1.375000,0.873580198) (1.437500,0.881239589) (1.500000,0.888434919) (1.562500,0.895194306) (1.625000,0.901544162) (1.687500,0.907509300) (1.750000,0.913113028) (1.812500,0.918377244) (1.875000,0.923322517) (1.937500,0.927968170) (2.000000,0.932332358) (2.062500,0.936432134) (2.125000,0.940283516) (2.187500,0.943901555) (2.250000,0.947300388) (2.312500,0.950493296) (2.375000,0.953492755) (2.437500,0.956310487) (2.500000,0.958957501) (2.562500,0.961444140) (2.625000,0.963780121) (2.687500,0.965974573) (2.750000,0.968036069) (2.812500,0.969972666) (2.875000,0.971791930) (2.937500,0.973500971) (3.000000,0.975106466) (3.062500,0.976614689) (3.125000,0.978031533) (3.187500,0.979362535) (3.250000,0.980612896) (3.312500,0.981787501) (3.375000,0.982890941) (3.437500,0.983927526) (3.500000,0.984901308) (3.562500,0.985816092) (3.625000,0.986675451) (3.687500,0.987482745) (3.750000,0.988241127) (3.812500,0.988953561) (3.875000,0.989622831) (3.937500,0.990251552) (4.000000,0.990842181) (4.062500,0.991397025) (4.125000,0.991918253) (4.187500,0.992407901) (4.250000,0.992867883) (4.312500,0.993299996) (4.375000,0.993705929) (4.437500,0.994087267) (4.500000,0.994445502) (4.562500,0.994782032) (4.625000,0.995098172) (4.687500,0.995395159) (4.750000,0.995674152) (4.812500,0.995936242) (4.875000,0.996182453) (4.937500,0.996413746) (5.000000,0.996631027)
};
\addlegendentry{$v_{1}$}
\addplot[thick,dashed] coordinates {
(0.000000,0.458333333) (0.062500,0.516426365) (0.125000,0.554649222) (0.187500,0.584575243) (0.250000,0.610599608) (0.312500,0.634371414) (0.375000,0.656514425) (0.437500,0.677281613) (0.500000,0.696797313) (0.562500,0.715143335) (0.625000,0.732387789) (0.687500,0.748593790) (0.750000,0.763821582) (0.812500,0.778128770) (0.875000,0.791570186) (0.937500,0.804197770) (1.000000,0.816060562) (1.062500,0.827204760) (1.125000,0.837673831) (1.187500,0.847508647) (1.250000,0.856747616) (1.312500,0.865426833) (1.375000,0.873580205) (1.437500,0.881239592) (1.500000,0.888434921) (1.562500,0.895194307) (1.625000,0.901544163) (1.687500,0.907509300) (1.750000,0.913113028) (1.812500,0.918377244) (1.875000,0.923322517) (1.937500,0.927968170) (2.000000,0.932332358) (2.062500,0.936432134) (2.125000,0.940283516) (2.187500,0.943901555) (2.250000,0.947300388) (2.312500,0.950493296) (2.375000,0.953492755) (2.437500,0.956310487) (2.500000,0.958957501) (2.562500,0.961444140) (2.625000,0.963780121) (2.687500,0.965974573) (2.750000,0.968036069) (2.812500,0.969972666) (2.875000,0.971791930) (2.937500,0.973500971) (3.000000,0.975106466) (3.062500,0.976614689) (3.125000,0.978031533) (3.187500,0.979362535) (3.250000,0.980612896) (3.312500,0.981787501) (3.375000,0.982890941) (3.437500,0.983927526) (3.500000,0.984901308) (3.562500,0.985816092) (3.625000,0.986675451) (3.687500,0.987482745) (3.750000,0.988241127) (3.812500,0.988953561) (3.875000,0.989622831) (3.937500,0.990251552) (4.000000,0.990842181) (4.062500,0.991397025) (4.125000,0.991918253) (4.187500,0.992407901) (4.250000,0.992867883) (4.312500,0.993299996) (4.375000,0.993705929) (4.437500,0.994087267) (4.500000,0.994445502) (4.562500,0.994782032) (4.625000,0.995098172) (4.687500,0.995395159) (4.750000,0.995674152) (4.812500,0.995936242) (4.875000,0.996182453) (4.937500,0.996413746) (5.000000,0.996631027)
};
\addlegendentry{$v_{6}$}
\addplot[thick,dashdotted] coordinates {
(0.000000,0.083333333) (0.062500,0.370673864) (0.125000,0.509373395) (0.187500,0.582710195) (0.250000,0.626893721) (0.312500,0.657686253) (0.375000,0.681970189) (0.437500,0.702771947) (0.500000,0.721442991) (0.562500,0.738606308) (0.625000,0.754566996) (0.687500,0.769490369) (0.750000,0.783479232) (0.812500,0.796607463) (0.875000,0.808934654) (0.937500,0.820512549) (1.000000,0.831387930) (1.062500,0.841603956) (1.125000,0.851200832) (1.187500,0.860216180) (1.250000,0.868685280) (1.312500,0.876641248) (1.375000,0.884115181) (1.437500,0.891136290) (1.500000,0.897732009) (1.562500,0.903928114) (1.625000,0.909748815) (1.687500,0.915216858) (1.750000,0.920353609) (1.812500,0.925179140) (1.875000,0.929712307) (1.937500,0.933970823) (2.000000,0.937971329) (2.062500,0.941729456) (2.125000,0.945259890) (2.187500,0.948576425) (2.250000,0.951692022) (2.312500,0.954618854) (2.375000,0.957368359) (2.437500,0.959951280) (2.500000,0.962377709) (2.562500,0.964657128) (2.625000,0.966798445) (2.687500,0.968810025) (2.750000,0.970699730) (2.812500,0.972474944) (2.875000,0.974142603) (2.937500,0.975709223) (3.000000,0.977180927) (3.062500,0.978563465) (3.125000,0.979862239) (3.187500,0.981082324) (3.250000,0.982228488) (3.312500,0.983305210) (3.375000,0.984316696) (3.437500,0.985266899) (3.500000,0.986159533) (3.562500,0.986998084) (3.625000,0.987785830) (3.687500,0.988525850) (3.750000,0.989221033) (3.812500,0.989874098) (3.875000,0.990487595) (3.937500,0.991063923) (4.000000,0.991605332) (4.062500,0.992113939) (4.125000,0.992591732) (4.187500,0.993040576) (4.250000,0.993462226) (4.312500,0.993858330) (4.375000,0.994230435) (4.437500,0.994579995) (4.500000,0.994908377) (4.562500,0.995216862) (4.625000,0.995506658) (4.687500,0.995778896) (4.750000,0.996034640) (4.812500,0.996274889) (4.875000,0.996500582) (4.937500,0.996712601) (5.000000,0.996911774)
};
\addlegendentry{$v_{13}$}
\addplot[thick,densely dotted] coordinates {
(0.000000,1.000000000) (0.062500,0.777442231) (0.125000,0.689875764) (0.187500,0.662807636) (0.250000,0.662436745) (0.312500,0.673458374) (0.375000,0.688969025) (0.437500,0.705877662) (0.500000,0.722821144) (0.562500,0.739217859) (0.625000,0.754838370) (0.687500,0.769610790) (0.750000,0.783532669) (0.812500,0.796631175) (0.875000,0.808945176) (0.937500,0.820517218) (1.000000,0.831390002) (1.062500,0.841604876) (1.125000,0.851201240) (1.187500,0.860216361) (1.250000,0.868685360) (1.312500,0.876641283) (1.375000,0.884115197) (1.437500,0.891136297) (1.500000,0.897732012) (1.562500,0.903928115) (1.625000,0.909748816) (1.687500,0.915216859) (1.750000,0.920353609) (1.812500,0.925179140) (1.875000,0.929712307) (1.937500,0.933970823) (2.000000,0.937971329) (2.062500,0.941729456) (2.125000,0.945259890) (2.187500,0.948576425) (2.250000,0.951692022) (2.312500,0.954618854) (2.375000,0.957368359) (2.437500,0.959951280) (2.500000,0.962377709) (2.562500,0.964657128) (2.625000,0.966798445) (2.687500,0.968810025) (2.750000,0.970699730) (2.812500,0.972474944) (2.875000,0.974142603) (2.937500,0.975709223) (3.000000,0.977180927) (3.062500,0.978563465) (3.125000,0.979862239) (3.187500,0.981082324) (3.250000,0.982228488) (3.312500,0.983305210) (3.375000,0.984316696) (3.437500,0.985266899) (3.500000,0.986159533) (3.562500,0.986998084) (3.625000,0.987785830) (3.687500,0.988525850) (3.750000,0.989221033) (3.812500,0.989874098) (3.875000,0.990487595) (3.937500,0.991063923) (4.000000,0.991605332) (4.062500,0.992113939) (4.125000,0.992591732) (4.187500,0.993040576) (4.250000,0.993462226) (4.312500,0.993858330) (4.375000,0.994230435) (4.437500,0.994579995) (4.500000,0.994908377) (4.562500,0.995216862) (4.625000,0.995506658) (4.687500,0.995778896) (4.750000,0.996034640) (4.812500,0.996274889) (4.875000,0.996500582) (4.937500,0.996712601) (5.000000,0.996911774)
};
\addlegendentry{$v_{24}$}
\end{axis}
\end{tikzpicture}
\caption{Four coordinates of the explicitly specified payoff-gradient game with initial state \eqref{eq:initial}. All 24 flows converge to one by Proposition~\ref{prp:gameequilibrium}. Curves are evaluations of \eqref{eq:gameSolution}; styles distinguish the selected players.}\label{fig:trajectories}
\end{figure}

\begin{table}[t]
\centering
\caption{Storage of the same interaction operator $A\tstp W$. The bridge is generated from its dimensions; damping and forcing are specified separately by the common scalars 13 and 1.}\label{tab:storage}
\begin{tabular}{@{}>{\ifrevisionblue\color{revisionblue}\fi}p{0.71\linewidth}>{\ifrevisionblue\color{revisionblue}\fi}r@{}}
\toprule
\textcolor{black}{Representation} & \textcolor{black}{Entries}\\
\midrule
\textcolor{black}{Rectangular tensor $A$} & \textcolor{black}{24}\\
\textcolor{black}{Square-slice t-product tensor $K$} & \textcolor{black}{36}\\
\textcolor{black}{Materialized block-circulant matrix $L$} & \textcolor{black}{144}\\
\textcolor{black}{Materialized dense matrix $I_2\otimes L$} & \textcolor{black}{576}\\
\textcolor{black}{Classical implementation retaining the factorization $A^{(i)}\Psi$} & \textcolor{black}{24}\\
\botrule
\end{tabular}
\end{table}

Table~\ref{tab:storage} compares stored entries, not different dynamical state dimensions. For an unrestricted $A$ of the same shape, the factorized family has 24 independent coefficients: $\Psi_{2\times3}$ has full row rank, so $A^{(i)}\mapsto A^{(i)}\Psi_{2\times3}$ is injective. An unrestricted square-slice tensor $K$ has 36 coefficients, but the equivalent $K$ here belongs to that 24-dimensional subspace. The materialized matrices do not acquire 144 or 576 independent parameters; they contain repetition and structural zeros. The particular numerical example fixes all of these coefficients.

The comparison establishes a smaller coefficient description within the chosen class. A classical implementation retaining the same factorization has the same storage economy. Thus no general runtime advantage over equally structured classical methods is claimed. For empirical use, the bridge and circular coupling should first be justified; coefficient estimation and approximation errors would then require data. The convergence result above belongs to the specified synthetic game and is not a general assertion.

\section{Concluding Remarks}\label{S9}

In this paper, we propose a framework based on DK-STP to generalize the t-product of cubic matrices and the corresponding dynamical systems.
To relax the inner-dimension compatibility requirement, we have extended the DK-STP to cubic matrices and described the resulting slice-wise control systems. To allow interactions between frontal-slice subsystems, we have introduced the t-STP, combining the circular coupling of the t-product with the bridge construction of the DK-STP. For a fixed number of frontal slices, we have established its algebraic properties and described the associated Lie algebra and Lie group. The t-STP contains the t-product as the case of matching inner dimensions and permits rectangular coefficient slices. A supply-network game illustrates the construction, its equivalent classical realization, and a globally convergent payoff-gradient adjustment law. The principal representational advantage is a smaller coefficient description for the specified operator class; the state dimension is unchanged. The construction therefore extends the coefficient description of t-product systems \cite{hoo21,mao24} while retaining their circular coupling, and complements transform-based tensor products \cite{ker15}, which vary the transform used to define the multiplication.

The resulting framework connects control systems over cubic matrices with structured classical systems. Its usefulness depends on whether the bridge factorization and circular slice coupling are appropriate for the application. Relative to the Einstein-product systems in \cite{cche21b}, the present formulation fixes a particular circular and bridge-factorized operator structure. The Lie structures express the standard commutator and unit-group constructions \cite{hal03} for this algebra, with the coefficient exponential linking them to the trajectories. Several questions remain for further investigation.

\begin{itemize}
\item[(i)] \textbf{Linear control problems.} The classical realization provides a route to studying controllability, observability, stability, and control design. A useful next step is to determine which results and design procedures preserve the prescribed tensor structure and admissible input channels, and when that structure can be exploited computationally. When $s=1$, the state space reduces to a matrix space, so the fixed-depth construction also applies to matrix-valued systems. Existing tensor-based system analyses \cite{hoo21,cche21b,mao24} provide reference formulations; the additional question here is whether a chosen design stays within the bridge-factorized coefficient family.
\item[(ii)] \textbf{Nonlinear control systems.} The space $\R^{m\times n\times s}$ is linearly isomorphic to $\R^{mns}$. Accordingly, geometric methods for nonlinear control systems \cite{isi95} can be applied when the maps defining the dynamics are well-defined smooth vector fields on this space. The product-generated vector fields considered here have additional algebraic structure. Understanding how this structure affects accessibility, feedback design, and computation is an open direction.
\item[(iii)] \textbf{Higher-order hypermatrices.} An array $A\in\R^{m_1\times\cdots\times m_d}$, with $d>3$, can be reshaped into a third-order array after grouping its indices. This suggests a way to extend the present constructions. However, the induced product and coupling depend on the chosen grouping, and their relation to the original multiway interactions requires further study. General tensor decompositions \cite{kol09} and products defined through other invertible transforms \cite{ker15} offer possible comparison frameworks, but compatibility with the chosen grouping would have to be established.
\item[(iv)] \textbf{Coefficient storage and computational performance.} For a state $X\in\R^{m\times q\times s}$, the rectangular coefficient tensor $A\in\R^{m\times p\times s}$ stores $mps$ entries. Its classical realization has additional block-circulant, Kronecker, and bridge structure, which should be retained in a fair computational comparison. Future work should quantify when this representation improves storage or runtime relative to equally structured alternatives, and assess the approximation error when the model is fitted to data. Extensions involving different numbers of frontal slices require a compatible lifting rule and a separate associativity analysis. Low-rank tensor numerical methods \cite{kho15} illustrate another route to economical representations. Comparing or combining them with the fixed bridge requires accounting for both approximation error and implementation cost.
\end{itemize}

\bmhead*{Acknowledgements}
This work is supported in part by the National Natural Science Foundation of China (NSFC) under Grant No. 62350037.

The authors are grateful to the anonymous reviewers for helpful suggestions and careful identification of errors.

\setlength{\bibsep}{0.8em}

\end{document}